\documentstyle{article}
\begin{document}

{\large \bf
\noindent
The Uniqueness of the Solution to Inverse Problems of 

\noindent
Interpolation of Positive Operators in Banach Lattices}
\footnote{The research was supported by the Russian Foundation for Basic 
Research, grant no.\ 98-01-00103.}

\bigskip
\noindent
O.\ E.\ TIKHONOV and L.\ V.\ VESELOVA

{\small
\noindent
Research Institute of Mathematics and Mechanics, 
Kazan State University, Universitetskaya Str.\ 17, Kazan, Tatarstan, 
420008 Russia

\noindent
E-mail: Oleg.Tikhonov@ksu.ru 

\medskip
\noindent
{\bf Abstract. }
We prove that an interpolation pair of Banach lattices is 
uniquely determined by the collection of intermediate spaces with 
the property that these
are interpolation spaces for positive operators. A correspondent 
result for exact interpolation is also presented. 


}

\bigskip
\noindent
In [2] we have proved that an arbitrary Banach couple (= 
interpolation pair of Banach spaces) is determined,
up to the equivalence of the norms, by the collection of all its 
interpolation spaces. By making use of this, we have proved 
in [3] that a Banach couple is determined by the collection of all 
its exact interpolation spaces up to the proportionality of the 
norms. (See also [4], where inverse interpolation problems were 
examined in a more general setting.) The goal of the present note 
is to show that similar assertions hold true when studying 
interpolation of positive operators in Banach lattices.

We first introduce notation and present some facts concerning
interpolation of linear operators in Banach lattices. For basic 
definitions of the theory we refer the reader to [1]. 

For a Banach lattice $E$, we denote by $\| \cdot \| _E$ and by $E^+$ 
the norm in $E$ and the cone of positive elements, respectively;
$|x|$ stands for the absolute value of $x \in E$. 
We denote by $\| \cdot \| ^E$ the norm in the dual Banach lattice 
$E^*$ and use the notation $\| T \| _{E \to E}$ to denote the norm 
of a bounded linear operator $T$ in $E$. A linear operator $T$ in 
$E$ is called {\it positive\/} if $T$ maps $E^+$ into $E^+$. 
Let $L^+ (E)$ stand for the cone of all such operators. Recall that a 
positive operator in a Banach lattice is automatically bounded. 

Let $E$ and $F$ be two Banach lattices with $E \hookrightarrow F$. 
(Here and subsequently, the notation $E \hookrightarrow F$ for two 
Banach lattices $E$ and $F$ means that $E$ is an ideal in $F$ and 
there exists a constant $c>0$ such that 
$\| x \| _E \le c \| x \| _F$ for all $x \in E$.) 
If $\varphi \in F^*$ then, clearly, $\varphi \big | _E \in E^*$. We 
write $\| \varphi \| ^E$ to denote $\| \varphi \big | _E \| ^E$. 
Note also that $| \varphi \big | _E | = | \varphi | \big | _E$.

Two Banach lattices $X$ and $Y$ are said to form {\it an 
interpolation pair $(X,Y)$ of Banach lattices\/} if they both are 
ideals in a certain vector lattice $U$ which is a Hausdorff 
topological vector space such that the embeddings of $X$ and of 
$Y$ into $U$ are continuous. To each interpolation pair $(X,Y)$ of 
Banach lattices can be canonically associated two Banach lattices, 
the intersection $X \cap Y$ and the sum $X+Y$. As can be easily 
seen,
$$
(X+Y)^+ = \{ x+y \mid x \in X^+ , y \in Y^+ \}. \eqno{(1)}
$$
We denote by $L^+ ((X,Y))$ the set of all linear operators $T$ in 
$X+Y$ with $T \big | _X \in L^+ (X)$, $T \big | _Y \in L^+ (Y)$.
One can easily deduce from (1) that $L^+ ((X,Y)) \subset L^+ (X+Y)$.
Also, we denote by $L_1 ((X,Y))$ the set of rank one linear 
operators $T$ in $X+Y$ with the property that $T$ maps $X$ into $X$ 
and $Y$ into $Y$ continuously. 

A Banach lattice $Z$ is called {\it intermediate\/} for an 
interpolation pair $(X,Y)$ of Banach lattices if 
$X \cap Y \hookrightarrow Z \hookrightarrow X+Y$. 
For a constant $c \ge 1$, we denote by $\hbox{P-Int}_c (X,Y)$ 
the collection of intermediate lattices $Z$ satisfying the following 
"interpolation property for positive operators": $T(Z) \subset Z$ and 
$$
\| T \| _{Z \to Z} \le 
c \max \{ \| T \| _{X \to X} , \| T \| _{Y \to Y} \} 
\eqno{(2)}
$$
for all $T \in L^+ ((X,Y))$. 
Set $\hbox{P-Int}(X,Y) = 
\bigcup \limits _{c \ge 1} \hbox{P-Int}_c (X,Y)$. 
Similarly, we denote by $\hbox{R1-Int}_c (X,Y)$ the collection of 
intermediate lattices $Z$ satisfying the following condition: 
$T(Z) \subset Z$ and (2) holds true for all 
$T \in L_1 ((X,Y))$. Set $\hbox{R1-Int}(X,Y) = 
\bigcup \limits _{c \ge 1} \hbox{R1-Int}_c (X,Y)$. 

\medskip
\noindent {\sc LEMMA 1.} 
{\it Let $(X,Y)$ be an interpolation pair of Banach lattices and 
$c \ge 1$. 
Then $\hbox{\rm P-Int}_c (X,Y) \subset \hbox{\rm R1-Int}_c (X,Y)$. 

\smallskip
Proof.} 
First, consider an intermediate Banach lattice $Z$ and a continuous linear 
operator $T$ in $X+Y$ of rank one. Necessary, $T$ is of the form 
$T( \cdot )= \varphi ( \cdot ) x$ with $\varphi \in (X+Y)^*$, 
$x \in X+Y$.
Define $|T|$ by $|T| ( \cdot ) = | \varphi | ( \cdot ) |x|$.
We claim that $|T|$ maps $Z$ into itself if and only if $T$ does, and 
$\| \, |T| \, \| _{Z \to Z} = \| T \| _{Z \to Z}$ 
when $T(Z) \subset Z$. 
Really, the following two cases are possible: $x \in Z$ and 
$x \notin Z$. In the first case we have 
$$ 
\| \, |T| \, \| _{Z \to Z} = 
\| \, | \varphi | \, \| ^Z \| \, |x| \, \| _Z =
\| \varphi \| ^Z \| x \| _Z = \| T \| _{Z \to Z} .
$$ 
In the second one,
$$ 
T(Z) \subset Z \quad \Longleftrightarrow \quad
\varphi \big | _Z = 0 \quad \Longleftrightarrow \quad  
| \varphi | \big | _Z = 0 \quad \Longleftrightarrow  \quad
|T| (Z) \subset Z 
$$
and  
$$
T(Z) \subset Z \quad \Rightarrow \quad \| T \| = \| \, |T| \, \| = 0 . 
$$
Thus, the claim is proved.

Now, let $Z \in \hbox{P-Int}_c (X,Y)$. Consider an operator 
$T \in L_1 ((X,Y))$. $T$ is clearly continuous in $X+Y$ and we can 
define $|T|$ as above. Then $|T| \in L^+ ((X,Y))$. Since 
$Z \in \hbox{P-Int}_c (X,Y)$, we have $|T| \big | _Z \in L^+ (Z)$. 
The latter implies that $T$ maps $Z$ into itself. Moreover, 
$$
\begin{array}{rl}
\| T \| _{Z \to Z} & = \| \, |T| \, \| _{Z \to Z} \\
& \le 
c \max \{ \| \, |T| \, \| _{X \to X} , \| \, |T| \, \| _{Y \to Y} \} 
\\  
& = c \max \{ \| T \| _{X \to X} , \| T \| _{Y \to Y} \} .
\end{array}
$$
Thus, $Z \in \hbox{R1-Int}_c (X,Y)$. 

\medskip
For two Banach lattices $E$ and $F$, let the notation $E \simeq F$ 
mean that $E$ and $F$ coincide as vector lattices and
their norms are equivalent, the notation $E\cong F$ 
mean that the norms are proportional.

For two interpolation pairs of Banach lattices $(X,Y)$ and $(V,W)$, 
we write $(X,Y) \simeq (V,W)$ if $X \simeq V$ and 
$Y \simeq W$ or $X \simeq W$ and $Y \simeq V$. 
We write $(X,Y) \cong (V,W)$ if $X \cong V$ and 
$Y \cong W$ or $X \cong W$ and $Y \cong V$. 

\medskip
\noindent {\sc THEOREM 1.} 
{\it 
Let $(X,Y)$ and $(V,W)$ be two interpolation pairs of Banach 
lattices with 
$$
\hbox{\rm P-Int}(X,Y) = \hbox{\rm P-Int}(V,W).
$$
Then $(X,Y) \simeq (V,W)$.

\smallskip
Proof.} 
Since $X,Y \in \hbox{P-Int}(X,Y)$ and $V,W \in \hbox{P-Int}(V,W)$, 
we obtain $X,Y \in \hbox{R1-Int}(V,W)$ and 
$V,W \in \hbox{R1-Int}(X,Y)$ by making use of Lemma 1. So it 
remains to prove 
$$
X,Y \in \hbox{R1-Int}(V,W) \hbox{ and } V,W \in \hbox{R1-Int}(X,Y) 
\quad \Rightarrow  \quad (X,Y) \simeq (V,W) .
$$
By analyzing proofs in [2] one can verify that this implications was 
actually proved there. Note also that from [4, Proposition 1.4 
and Theorem 2.8] the latter implication follows immediately. 
  
\medskip
\noindent {\sc THEOREM 2.} 
{\it 
Let $(X,Y)$ and $(V,W)$ be two interpolation pairs of Banach 
lattices with 
$$
\hbox{\rm P-Int} _1 (X,Y) = \hbox{\rm P-Int} _1 (V,W).
$$
Then $(X,Y) \cong (V,W)$.

\smallskip
Proof.} 
We can proceed similarly to the proof of Theorem 1, but here we need 
to refer to [3] and to [4, Theorem 3.7]. 

\bigskip
\noindent
{\bf References}

\smallskip
1. J.~Linderstrauss and L.~Tzafriri, {\it Classical Banach Spaces II. 
Function Spaces,} Springer-Verlag, Berlin-Heidelberg-New York, 1979.

2. O.~E.~Tikhonov and L.~V.~Veselova, A Banach couple is determined 
by the collection of its interpolation spaces, {\it Proc.\ Amer.\ 
Math.\ Soc.} {\bf 126} (1998), 1049--1054. 

3. O.~E.~Tikhonov and L.~V.~Veselova, The uniqueness of the solution 
to the inverse problem of exact interpolation, {\it Israel Math.\ 
Conf.\ Proc.} {\bf 13} (1999), 209--215.

4. L.~V.~Veselova and O.~E.~Tikhonov, {\it The uniqueness of the 
solution to inverse interpolation problems}, Research Institute of 
Mathematics and Mechanics, Preprint no.\ 95-2, Kazan Mathematics 
Foundation, Kazan, 1995 (Russian); English transl., XXX e-Print Archive,  
math.FA/9902108.  
 
\end{document}